# Аналитическое конструирование регуляторов для колебательных систем с жидкими демпферами

**Алиев[1] Ф.А., Алиев[1] Н.А., Исмаилов[1] Н.А.**
[1]Институт Прикладной Математики, БГУ, Баку, Азербайджан

**Абстракт.** Рассматривается математическая модель управления колебательных систем с жидкими демпферами, которая отличается от классических заменой первой производной такой производной дробного порядка, которая находится между числами 0 и 2, отличными от единицы. Используя методы построения регуляторов Летова строится закон регулирования, обеспечивающий асимптотическую устойчивость замкнутой системы и минимизирующий квадратичный функционал. Результаты иллюстрируются конкретным числовым примером.

**Введение.** После появления классической работы Летова [1] многие современные задачи из техники- задачи виброзащиты [2,3], робототехники [4-8], ядерные реакторы [9,10] и др. на основе этого метода бурно развивались. В этих задачах движение объекта описывается классическими обыкновенными дифференциальными уравнениями с натуральными числами порядка производных. Однако существуют многие задачи из нефтяной индустрии, металловедения, колебательной системы, где порядки производных являются дробными числами, которые требует отдельный подход [11-13], где нахождение регуляторов, стабилизирующие эти процессы, требует отдельных разработок [14].

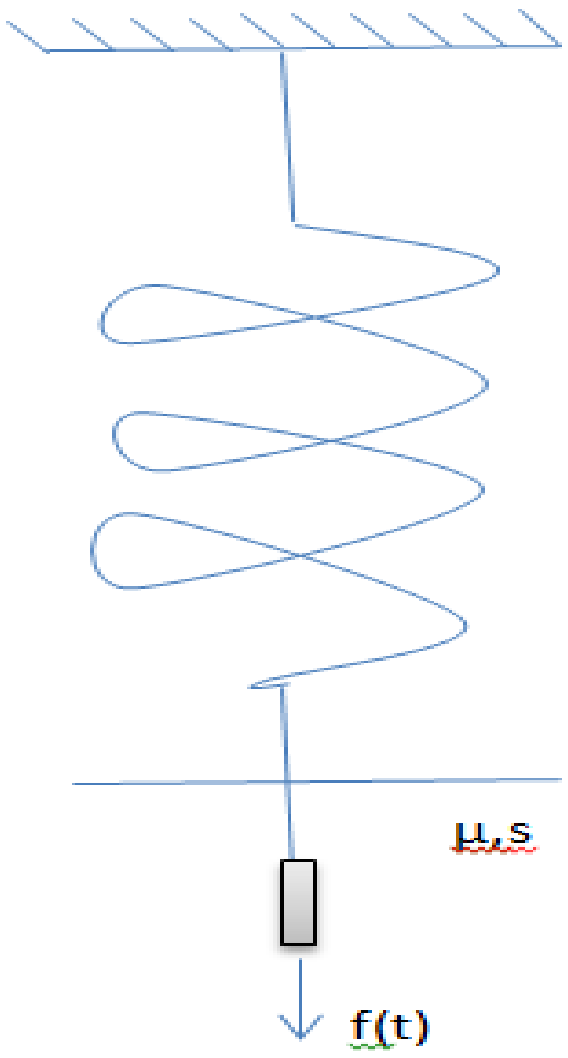

фиг.1

В данной работе, на основе [1,14], приводится конкретный вид регуляторов, обеспечивающих стабилизацию колебательных систем с жидкими демпферами и минимизирующих соответствующий квадратичный функционал. Для получения конкретного вида уравнения регуляторов, исходное дифференциальное уравнение с дробной производной сводится (с шагом, равном знаменателю) к нормальной системе дифференциальных уравнений дробного порядка. Результаты иллюстрируются примером, где дробный порядок является рациональным числом.

**2. Постановка задачи.** Пусть движение колебательных систем с жидкими демпферами [13-15] описывается (см. Фиг.1) следующими обыкновенными дифференциальными уравнениями с дробной производной

$$y''(x)+aD^{\alpha}y+by=u \qquad (0<\alpha<2,\ \alpha\neq 1) \tag{1}$$

с начальными условиями

$$y(x_0)=0,\quad y'(x_0)=y_1, \tag{2}$$

где $a=\dfrac{2s\sqrt{\mu\rho}}{m}$, $b=\dfrac{k}{m}$, и рассматривается жесткая пластина с массой $m$ и площадью $s$, $\rho$-плотность жидкости, $\mu$-постоянная вязкость упругости, постоянная $k$ характеризует свойства пружины (см. Фиг.1), $u$ - управляющее воздействие, являющееся внешней силой.

Задача состоит в нахождении такого закона регулирования

$$u=Wy\ , \tag{3}$$

чтобы замкнутая система (1)-(3) была асимптотически устойчивой и квадратичный функционал

$$J=\frac{1}{2}\int_0^{\infty}\left(qy^2(t)+ru^2(t)\right)dt \tag{4}$$

получил минимальное значение. Для классической задачи синтеза (1)-(4) существуют два метода-временной метод [1, 16-18] и частотной метод [19-25], которые требуется модифицировать для данного случая. Поскольку этой модификацией мы можем выявить структуру W из (3), в данной работе остановимся на создании временного метода, т.е. аналога АКОР Летова А.М. для стабилизации колебательных систем с жидкими демпферами около соответствующей программной траектории и управлений.

**3. Приведение задачи (1)-(4) к задачам для нормальной системы**. Для того, чтобы использовать результаты [1, 11, 14] приведем уравнение (1) к нормальному виду. Поэтому, используя результаты [13], приведем операторное уравнение (1) в общем случае к нормальной системе уравнений. Для простоты предположим, что $\alpha$ является несократимым рациональным числом $\alpha=\dfrac{p}{q}$, где $p$ и $q$ нечетные натуральные числа и $p<2q$, иначе $\alpha>2$. Если для производной уравнения (1) принять шаг, равный $\dfrac{1}{q}$, то уравнение (1) примет вид

$$(D^{\frac{1}{q}}y)^{2q}+a\left(D^{\frac{1}{q}}\right)^{p}y+by=u. \tag{5}$$

Теперь примем следующее обозначение

$$y = z_1,\ D^{\frac{1}{q}} y = D^{\frac{1}{q}} z_1 = z_2,\quad D^{\frac{2}{q}} y = D^{\frac{1}{q}} z_2 = z_3, \dots, D^{\frac{p-1}{q}} y = D^{\frac{1}{q}} z_{p-1} = z_p, \dots,$$

$$D^{\frac{p}{q}} y = D^{\frac{1}{q}} z_p = z_{p+1,\dots,}\ y' = z_1'' = \left(D^{\frac{1}{q}}\right)^{2q} z_1 = D^{\frac{1}{q}} z_{2q} \qquad (6)$$

и, учитывая (6) в соотношении (3), получаем последнее (5) в следующем нормальном виде

$$D^{\frac{1}{q}} \begin{bmatrix} z_1 \\ z_2 \\ z_3 \\ \dots \\ \dots \\ z_{p-1} \\ z_p \\ z_{p+1} \\ \dots \\ \dots \\ z_{2q} \end{bmatrix} = \begin{bmatrix} 0 & 1 & 0 & 0 \dots\dots\dots\dots 0 & 0 & 0 \\ 0 & 0. & 1 & 0 \dots\dots\dots\dots 0 & 0 & 0. \\ .0 & 0 & 0 & 1 \dots\dots\dots\dots 0 & 0 & 0 \\ \dots\dots & & & & & \\ \dots\dots & & & & & \\ \dots\dots & & & & & \\ \dots\dots & & & & & \\ \dots\dots & & & & & \\ 0 & 0 & 0 & 0 \dots\dots\dots\dots 1 & 0 & 0 \\ 0 & 0 & 0 & 0 \dots\dots\dots\dots 0 & 1 & 0 \\ 0 & 0 & 0 & 0 \dots\dots\dots\dots 0 & 0 & 1 \\ \underbrace{-b\ 0\ 0\ 0 \quad -a}_{p} & & & & 0 & 0 & 0 \end{bmatrix} \begin{bmatrix} z_1 \\ z_2 \\ z_3 \\ \dots \\ \dots \\ z_p \\ z_{p+1} \\ \dots \\ \dots \\ \dots \\ z_{2q-1} \\ z_{2q} \end{bmatrix} + \begin{bmatrix} 0 \\ 0 \\ 0 \\ \dots \\ \dots \\ 0 \\ 0 \\ \dots \\ \dots \\ 0 \\ 0 \\ 1 \end{bmatrix} u = Fz + Gu \qquad (7)$$

т.е. (7) и (5) или (1) эквивалентны. Таким образом, как в [14], можно ставить следующую задачу аналитического конструирования регуляторов Летова А.М. [16] для колебательных систем (5) в следующем виде; требуется найти закон регулирования

$$u = k\bar{z} \qquad (8)$$

такой, чтобы минимизировался критерий качества (2), являющийся эквивалентным

$$J = \frac{1}{2}\int_0^\infty \left(z'Qz + ru^2\right)dt, \qquad (9)$$

при котором замкнутая система (7)-(8) являлась бы асимптотически устойчивой, где $Q = Q' \geq 0,$ - известная матрица, штрих означает операцию транспонирования. Отметим, что здесь $Q = \left(\left[q_{i,y}\right]_{i,j=1}^{2q}\right),\quad q_{11} = q,\ \ q_{i,j} = 0.$

**4. Временной метод.** Используя результаты [14], можем найти закон регулирования (8) в следующем виде

$$u = -\frac{1}{r} G' T_3 T_1^{-1} z, \qquad (10)$$

где $T_1, T_3$ определяются в следующем виде

$$T^{-1} H T = \begin{bmatrix} A_+ & 0 \\ 0 & A_- \end{bmatrix}, \quad T = \begin{bmatrix} T_1 & T_2 \\ T_3 & T_4 \end{bmatrix}, \quad H = \begin{bmatrix} F & -\frac{1}{r} G G' \\ -Q & -F' \end{bmatrix}, \qquad (11)$$

где $\lambda(A_+) < 0,\ \lambda(A_-) > 0,$ являются собственными значениями матрицы $H$ ,а $H$ - Гамильтановая матрица.

Используя результаты [1, 14, 16, 17], легко показать, что замкнутая система (7)-(10)

$$D^{\frac{1}{q}} z = (F - \frac{1}{r} G' T_3 T_1^{-1}) z, \quad z(0) = z_0 \qquad (12)$$

является асимптотически устойчивой. Обозначая

$$T_3 T_1^{-1} = \left[ t_{ij} \right]_{i,j'=1}^{2q} \qquad (13)$$

для регулятора (10) имеем следующий вид

$$u = -\frac{1}{r} \left[ t_{2q,1}\ y + t_{2q,2} D^{\frac{1}{q}} y + t_{2q,3} D^{\frac{2}{q}} y + \ldots. + t_{2q,p+1} D^{\frac{p}{q}} y + \ldots + t_{2q,2q} D^{\left(2 - \frac{1}{q}\right)} y \right]. \qquad (14)$$

Отметим, что в (13) $T_3 T_1^{-1}$ является положительно определенным решением следующего матричного алгебраического уравнения Риккати (АУР) [18, 26, 27]

$$SF + F'S - SGC^{-1}G'S + R = 0, \qquad (15')$$

где существует много эффективных вычислительных алгоритмов для ее решения [18, 27, 29, 30].

Таким образом, вследствие свойства замкнутой системы (12) аналогичное уравнение (5)-(14) имеет вид

$$y'' + t_{2q,2q} D^{\left(2 - \frac{1}{q}\right)} y + \ldots\ldots + \left( a + \frac{1}{r} t_{2q,p+1} \right) D^{\frac{p}{q}} y + \ldots. + \frac{1}{r} t_{2q,2} D^{\frac{1}{q}} y + \left( b + \frac{1}{r} t_{2q,1} \right) y = 0, \qquad (16)$$

которое является асимптотически устойчивой.

Отметим, что одни из параметров $p$ и $q$, если они являются четными, то ситуация осложняется и требуется отдельное исследование. Действительно, пусть $p = 2k$, а $q = 2m + 1$. Тогда $\alpha = \frac{2k}{2m+1}$ преобразуем следующим образом

$$\alpha = \frac{2k(2l+1)}{(2m+1)(2l+1)} \approx \frac{4kl + 2k + 1}{4ml + 2(m+l) + 1}.$$

Таким образом $\alpha$ приближенно представлено в виде отношения двух нечетных чисел, которое было рассмотрено выше. Приведем следующий пример, который иллюстрирует вышеприведенные результаты.

**5.1. Пример.** Пусть в (1) $\alpha=\frac{1}{3},\ a=3,\ b=1$, а в функционале (4) $q=1,\ r=1$. Тогда легко можем уравнение (3) в данном случае написать в виде

$$(D^{\frac{1}{3}})^6 y + D^{\frac{1}{3}} y + y = u. \tag{17}$$

И нормальная система с шагом 1/3 для (14) будет

$$D^{\frac{1}{3}} z = Fz + Gu, \tag{18}$$

где

$$z=\begin{bmatrix} z_1 \\ z_2 \\ \cdots \\ \cdots \\ \cdots \\ z_6 \end{bmatrix}, \qquad F=\begin{bmatrix} 0_{5\times1} & E_{5\times5} \\ -1 & f_1 \end{bmatrix}, \qquad f_1'=\begin{bmatrix} -3 \\ 0 \\ 0 \\ 0 \\ 0 \end{bmatrix},$$

$E_{5\times5}$ - единичная матрица порядка $5\times5$, $0_{5\times1}$ - нулевой вектор, $G'=[0,0,0,0,0,1]$.

Тогда $Q=\begin{bmatrix} 1 & 0_{1\times5} \\ 0_{5\times1} & 0_{5\times5} \end{bmatrix}$. Составляем Гамильтоновы матрицы из (11) и, используя пакет программ MATLAB для $T_1,\ T_3$, имеем конкретные числовые матрицы, а $T_3\cdot T_1^{-1}$ имеет следующий вид

$$T_3\cdot T_1^{-1}=S, \tag{19}$$

где $S$ является положительно определенным решением матричного АУР [18], а уравнение регулятора составляется через элементы $t_{6,j}\ \left(j=\overline{1,6}\right)$ из (19) в виде

$$\begin{aligned} u=&-0.4142y-3.5878D^{\frac{1}{3}}y-12.162D^{\frac{2}{3}}y \\ &-13.2864Dy-9.5964D^{\frac{4}{3}}y-4.3809D^{\frac{5}{3}}y. \end{aligned} \tag{20}$$

Учитывая (20) в (17) имеем следующую замкнутую систему (17)-(20)

$$\begin{aligned} &y''+4.3809D^{\frac{5}{3}}y+9.5964D^{\frac{4}{3}}y+13.2864Dy+12.162D^{\frac{2}{3}}y \\ &+6.5878D^{\frac{1}{3}}y+1.4142y=0. \end{aligned} \tag{21}$$

Для доказательства асимптотической устойчивости решения (21), исходя из функции Миттага-Леффлера [11,31], решение операторного дифференциального уравнения (21) с шагом 1/3 представим в следующем виде

$$y(x,\lambda)=\sum_{k=-2}^{\infty}\lambda^{k+2}\frac{x^{\frac{k}{3}}}{\left(\frac{k}{3}\right)!} \qquad (22)$$

Легко видеть, что при $x>0$, подставляя соотношение (23)

$$D^{\frac{m}{3}}y=\lambda y,(m=\overline{1,\ 6}). \qquad (23)$$

и (22) в (21) , получим следующее характеристическое уравнение

$$\lambda^6+4.3809\lambda^5+9.5964\lambda^4+13.2864\lambda^3+12.162\lambda^2+6.5878\lambda+1.4142=0. \qquad (24)$$

С помощью MATLAB находим $\lambda_i$ $(i=\overline{1,6})$, для которого $\operatorname{Re}\lambda_i<0$, т.е. замкнутая система (21) асимптотически устойчива.

Отметим, что общее решение уравнения (21) из [31, 11] имеет вид

$$y(x)=\sum_{l=1}^{6}c_l\sum_{k=-2}^{\infty}\lambda_l^{k+2}\frac{x^{\frac{k}{3}}}{\left(\frac{k}{3}\right)!}. \qquad (25)$$

Начальные условия для (25) определяются в следующем виде

$$y(x_0)=0,\ D^{\frac{1}{3}}y(x_0)=0,\ D^{\frac{2}{3}}y(x_0)=0,\ y'(x_0)=y_1,\ D^{\frac{4}{3}}y(x_0)=0,\ D^{\frac{5}{3}}y(x_0)=0. \qquad (26)$$

а $c_l\ (l=\overline{1,\ 6})$ определяется в следующем виде [11, 31, 32]

$$c_s=\frac{\Delta_S}{\Delta},\quad s=\overline{1,\ 6}, \qquad (27)$$

где

$$\Delta=\left(\prod_{i=1}^{6}y(x_0,\lambda_i)\right)\begin{vmatrix}1&1&1&1&1&1\\ \lambda_1&\lambda_2&\lambda_3&\lambda_4&\lambda_5&\lambda_6\\ \lambda_1^2&\lambda_2^2&\lambda_3^2&\lambda_4^2&\lambda_5^2&\lambda_6^2\\ \dots&\dots&\dots&\dots&\dots&\dots\\ \dots&\dots&\dots&\dots&\dots&\dots\\ \dots&\dots&\dots&\dots&\dots&\dots\\ \lambda_1^5&\lambda_2^5&\lambda_3^5&\lambda_4^5&\lambda_5^5&\lambda_6^5\end{vmatrix},\ \Delta_S=\left(\prod_{\substack{i=1\\ i\neq S}}^{6}y(x_0,\lambda_i)\right)\begin{vmatrix}1&1&\dots&0&\dots&1\\ \lambda_1&\lambda_2&\dots&0&\dots&\lambda_6\\ \lambda_1^2&\lambda_2^2&\dots&0&\dots&\lambda_6^2\\ \dots&\dots&\dots&y_1&\dots&\dots\\ \dots&\dots&\dots&0&\dots&\dots\\ \dots&\dots&\dots&\dots&\dots&\dots\\ \lambda_1^5&\lambda_2^5&\dots&0&\dots&\lambda_6^5\end{vmatrix}.$$

Тогда окончательный вид решения замкнутой системы (21) с начальными условиями (26) примет вид

$$y(x)=\sum_{k=1}^{6}\frac{\Delta_k}{\Delta}\left[\sum_{s=0}^{1}\lambda_k^s\frac{d}{dx}\int_0^x\frac{(x-t)^{\frac{s-2}{3}}}{\frac{s-2}{3}!}e^{\lambda_k^3 t}dt+\lambda_k^2 e^{\lambda_k^3 x}\right]. \qquad (28)$$

Поскольку в (28) подынтегральное выражение имеет слабую особенность, то производные нельзя перевести под знаком интеграла. Поэтому это выражение преобразуем следующим образом

$$\frac{d}{dx}\int_0^x\frac{(x-t)^{\frac{s-2}{3}}}{\frac{s-2}{3}!}e^{\lambda_k^3 t}dt=-\frac{d}{dx}\int_0^x e^{\lambda_k^3 t}d_t\frac{(x-t)^{\frac{s+1}{3}}}{\frac{s+1}{3}!}=-\frac{d}{dx}\left[e^{\lambda_k^3 t}\frac{(x-t)^{\frac{s+1}{3}}}{\frac{s+1}{3}!}\Bigg|_{t=0}^{x}-\int_0^x\lambda_k^3 e^{\lambda_k^3 t}\frac{(x-t)^{\frac{s+1}{3}}}{\frac{s+1}{3}!}dt\right]$$

$$=\frac{x^{\frac{s-2}{3}}}{\frac{s-2}{3}!}+\lambda_k^3\int_0^x e^{\lambda_k^3 t}\frac{(x-t)^{\frac{s-2}{3}}}{\frac{s-2}{3}!}dt. \qquad (29)$$

Подставляя (29) в (28) имеем

$$y(x)=\sum_{k=1}^{6}\frac{\Delta_k}{\Delta}\left[\sum_{s=0}^{1}\lambda_k^s\left[\frac{x^{\frac{s-2}{3}}}{\frac{s-2}{3!}}+\lambda_k^3\int_0^x e^{\lambda_k^3 t}\frac{(x-t)^{\frac{s-2}{3}}}{\frac{s-2}{3}!}dt\right]+\lambda_k^2 e^{\lambda_k^3 x}\right] \qquad (30)$$

решение замкнутой системы (5)-(14) в данном случае уравнения (21).

Из-за $\mathrm{Re}\,\lambda_i<0$ $(i=\overline{1,\ 6})$ $y(x)$ из (30) стремится к нулю при $x\to\infty$. Действительно, первое слагаемое в (30) содержит отрицательную степень $x$. Точно такое же слагаемое содержится под знаком интеграла. А последнее слагаемое экспоненциальной функции с отрицательным показателем. Поэтому при $x\to\infty$ решение (30) $y(x)\to 0$, т.е. замкнутая система (21) асимптотически устойчива.

**5.2.** В нашем предыдущем примере был рассмотрен случай, когда $\alpha=\frac{p}{q}$, где $p$ и $q$- несократимые нечетные числа. Как доказано в [11, 14, 33] результаты построения регуляторов [1] только в таком случае верны. Когда одно из $p$ или $q$ является четным, то свойства зеркальности собственных значений матрицы $H$ и Калмана [1, 16] из (11) теряют смысл и появляется ощущение, что теория Летова теряет смысл. Однако это не так. Как в [11, 32] показано, что для любых четных $\frac{p}{q}$ можно ее сколь угодно точно аппроксимировать нечетным $\frac{\tilde p}{\tilde q}$. Таким образом, можно решать задачи АКОР Летова [1] для любых случаев $p$ и $q$. Проиллюстрируем это на следующем простом примере, т.е. пусть в (1) $m=0,\ a=1,\ b=-1,\ \alpha=\frac{1}{2}$.
Тогда задачиу АКОР в данном случае можно сформулировать в виде:
задается динамическая система

$$D^{\frac{1}{2}}y=y+u,\ y(t_0)=1 \tag{31}$$

с квадратичным функционалом

$$J=\frac{1}{2}\int_{t_0}^{\infty}\left(3y^2+u^2\right)dt. \tag{32}$$

Если составим Гамильтонову матрицу АКОР (31), (32) то ее собственные значения будут $\mu_{1,2}=\pm 2i$, а это не обеспечивает асимптотическую устойчивость соответствующей замкнутой системы. Поэтому решаем задачи (31), (32) для $\tilde{\alpha}=\frac{2k+1}{2(2k+1)+1}\approx\frac{1}{2}$ при достаточно больших значениях $\kappa$.
Теперь имеем вместо уравнения (31) следующее уравнение

$$D^{\frac{2k+1}{4k+3}}y=y+u. \tag{33}$$

В таком случае для задачи (33), (32) собственные значения матрицы $H$ имеют вид $\tilde{\mu}_{1,2}=\pm 2,$ а собственные значения соответствующих гамильтоновых уравнений задачи (33), (32) будут

$$\lambda_1=2^{\frac{4k+3}{2k+1}},\quad \lambda_2=\left(-2\right)^{\frac{4k+3}{2k+1}},$$

т.е. зеркальность сохраняется. Легко представить, что в данном случае $T$ из (11) будет

$$T=\begin{bmatrix}1 & 1\\ 3 & -1\end{bmatrix},\ T^{-1}=\begin{bmatrix}\frac{1}{4} & \frac{1}{4}\\ \frac{3}{4} & -\frac{1}{4}\end{bmatrix}, \tag{34}$$

а уравнение регулятора (10) имеет вид

$$u(t)=-3y(t), \tag{35}$$

который обеспечивает асимптотическую устойчивость замкнутой системы

$$D^{\frac{2k+1}{4k+3}}y=-2y(t) \tag{36}$$

Теперь вычислим решение (32) в следующем виде [11, 32, 34]

$$y(t)=\left\{\left[\sum_{s=0}^{4k+2}(-2)^{\frac{s+4k+3}{2k+1}}\int_0^{t_0}\frac{\left(t_0-\tau\right)^{\frac{s-4k-2}{4k+3}}}{\frac{s-4k-2}{4k+3}!}e^{\tau^{(-2)\frac{4k+3}{2k+1}}}d\tau+\sum_{s=0}^{4k+2}(-2)^{\frac{s}{2k+1}}\frac{t_0^{\frac{s-4k-2}{4k+3}}}{\frac{s-4k-2}{4k+3}!}\right]^{-1}\right.$$

$$\left.\times\left[\sum_{s=0}^{4k+2}(-2)^{\frac{s+4k+3}{2k+1}}\int_0^{t}\frac{\left(t-\tau\right)^{\frac{s-4k-2}{4k+3}}}{\frac{s-4k-2}{4k+3}!}e^{\tau^{(-2)\frac{4k+3}{2k+1}}}d\tau+\sum_{s=0}^{4k+2}(-2)^{\frac{s}{2k+1}}\frac{t^{\frac{s-4k-2}{4k+3}}}{\frac{s-4k-2}{4k+3}!}\right]\right\}, \tag{37}$$

где там же показано, что при $t \to \infty$ решение (37) $y(t) \to 0$. Действительно, поскольку первый множитель в (37) не зависит от $t$, преобразуем второй множитель следующим образом:

$$[\bullet]=\left[\sum_{s=0}^{4k+1}(-2)^{\frac{s+4k+3}{2k+1}}\int_0^t \frac{(t-\tau)^{\frac{s-4k-2}{4k+3}}}{\frac{s-4k-2}{4k+3}!}e^{\tau^{(-2)\frac{4k+3}{2k+1}}}d\tau+(-2)^{\frac{8k+5}{2k+1}}\int_0^t e^{\tau^{(-2)\frac{4k+3}{2k+1}}}d\tau+\sum_{s=0}^{4k+1}(-2)^{\frac{s}{2k+1}}\frac{t^{\frac{s-4k-2}{4k+3}}}{\frac{s-4k-2}{4k+3}!}+4\right],(38)$$

где $[\bullet]$ означает второй множитель в правой части (37). Первое слагаемое в правой части (38) из-за отрицательности степени подынтегральной функции $(t-\tau)$ при $t \to \infty$ стремится к нулю. Точно так же предпоследнее слагаемое в (38) тоже при $t \to \infty$ стремится к нулю. Теперь вычислим интеграл, входящий во второе слагаемое правой части соотношения (38)

$$(-2)^{\frac{8k+5}{2k+1}}\int_0^t e^{\tau^{(-2)\frac{4k+3}{2k+1}}}d\tau=(-2)^{\frac{8k+5}{2k+1}}\left.\frac{e^{\tau^{(-2)\frac{4k+3}{2k+1}}}}{(-2)^{\frac{4k+3}{2k+1}}}\right|_{\tau=0}^{t}=(-2)^{\frac{4k+2}{2k+1}}\left[e^{t^{(-2)\frac{4k+3}{2k+1}}}-1\right]=4e^{t^{(-2)\frac{4k+3}{2k+1}}}-4.$$

Полученное первое слагаемое при $t \to \infty$ стремится к нулю. А второе слагаемое сокращается с последним слагаемым в правой части соотношения (38). То есть замкнутая система (36) асимптотически устойчива.

Теперь проведем вычисление при $k=1,\,2$ решения (37), которое показано на следующем рисунке.

Таким образом, при $k \to \infty$, $\alpha \to \frac{1}{2}$ и решение $y(t)$ из (37) стремится к $\tilde{y}(t)=e^{-4(t-t_0)}y(t_0)$, которое не является решением замкнутой системы (36) при $k \to \infty$, т.е. (37) можно принять за приближенное решение при фиксированном большом к замкнутой системы (36), а оптимальным регулятором для задачи АКОР (31), (32) является (35).

**Заключение.** Решается задача аналитического конструирования регуляторов Летова А.М. для колебательных систем с жидкими демпферами. Показывается, что соответствующая замкнутая система является асимптотически устойчивой. Результаты иллюстрируются числовыми примерами.

## Литература

## Analytical construction of regulators for oscillatory systems with liquid dampers

**Aliev F.A., Aliev N.A., Ismailov N.A.**

**Abstract.** A mathematical model of oscillatory systems control with liquid dampers is considered which differs from the classical ones by replacing the first derivative of such a fractional derivative that is between numbers 0 and 2 other than unity. Using the methods of constructing Letov controllers, a control law is constructed that ensures the asymptotic stability of the closed system and minimizes the quadratic functional. The results are illustrated by a specific numerical example.